\RequirePackage{fix-cm}
\documentclass[smallextended]{svjour3_MultibandJamming}       
\smartqed  
\usepackage{graphicx}
\usepackage{amsmath}
\usepackage{amssymb}
\usepackage{url}
\usepackage{algpseudocode}
\usepackage{algorithm}
\journalname{Optimization Letters}

\begin{document}

\title{
Revisiting wireless network jamming by SIR-based considerations and Multiband Robust Optimization*
\thanks{*This is the author's final version of the paper published in Optimization Letters 9(8), 1495-1510, 2015, DOI: 10.1007/s11590-014-0839-2r \; .
The final publication is available at Springer  http://dx.doi.org/10.1007/s11590-014-0839-2.}
}

\titlerunning{Multiband Robust Optimization for Wireless Network Jamming}        

\author{
Fabio D'Andreagiovanni
}


\institute{Fabio D'Andreagiovanni \at
            - Dept. of Optimization, Zuse Institute Berlin (ZIB), Takustr. 7, 14195 Berlin, Germany\\
            - DFG Research Center MATHEON, Technical University Berlin,
            Stra{\ss}e des 17. Juni 135, 10623 Berlin, Germany\\
            - Einstein Center for Mathematics,
            Stra{\ss}e des 17. Juni 135, 10623 Berlin, Germany\\
              Tel.: +49 30 84185460, Fax: +49 30 84185269, \email{d.andreagiovanni@zib.de}
}

\date{Received: 11 April 2014 / Accepted: 29 November 2014}

\maketitle

\begin{abstract}
We revisit the mathematical models for wireless network jamming introduced by Commander et al. \cite{CoEtAl07,CoEtAl08}: we first point out the strong connections with classical wireless network design and then we propose a new model based on the explicit use of signal-to-interference quantities. Moreover, to address the uncertain nature of the jamming problem and tackle the peculiar right-hand-side (RHS) uncertainty of the corresponding model, we propose an original robust cutting-plane algorithm drawing inspiration from Multiband Robust Optimization. Finally, we assess the performance of the proposed cutting plane algorithm by experiments on realistic network instances.
\end{abstract}

\section{Introduction}
\label{intro}

Wireless network jamming consists  in compromising the functionality of a wireless network by activating jamming devices (\emph{jammers}) that disrupt network communications by emitting interfering signals on the same frequencies of the network. Jamming is commonly associated with military and security questions: it is immediate to think about jamming hostile networks in war scenarios, to deny enemy communications, or in high-risk events, such as visits of heads of State, during which it is imperative to avoid bomb detonation by cellular phones. However, in recent times, jamming is also increasingly found in peaceful everyday-life applications that are not related to military and security issues. Italy provides two examples of such applications: the Italian public administration has evaluated the expediency of using jammers during large competitive examinations to prevent cheating, while schools have started to deploy jammers to avoid that students get distracted by smart phones during lectures\footnote{
Panorama. Cellulari a scuola: la soluzione c'\`e ma la vietano (in Italian).
\url{http://italia.panorama.it/Cellulari-a-scuola-la-soluzione-c-e-ma-la-vietano}
(2007)}.
USA provides another example: in some hotels, there is the suspicion that unscrupulous managers have shrewdly placed jammers to deny cellular coverage and force businessmen to use room phones, in an attempt to raise the final bill of stays\footnote{
C. Elliott: The Cellphone That Doesn't Work at the Hotel.
\url{http://www.nytimes.com/2004/09/07/business/07jamming.html?_r=0}.
The New York Times 07.09.2004
(2004)}.

The \emph{Wireless Network Jamming Problem (NJP)} can be described as the problem of optimally placing and configuring a set of jammers in order to interdict communications of a wireless network.
As pointed out by Commander et al. in \cite{CoEtAl07,CoEtAl08}, though the problem is very relevant and there is a wide literature about \emph{preventing} hostile jamming, surprisingly the NJP has been practically neglected until their work. Moreover, together with the work \cite{CoEtAl09}, these seem to be the only papers that have directly addressed the problem.

\noindent
Our main original contributions in the present paper are:
\begin{enumerate}
\item revisiting the models for the NJP introduced by Commander et al. \cite{CoEtAl07,CoEtAl08,CoEtAl09}. Specifically, we highlight the strong connections of the NJP with classical wireless network design and, as recommended by regulatory bodies, we adopt a testpoint model and signal-to-interference (SIR) quantities
to represent coverage and jamming conditions, refining the models of \cite{CoEtAl07,CoEtAl08,CoEtAl09};
\item addressing the uncertain nature of the NJP, considering a more realistic in-between case w.r.t. \cite{CoEtAl07} (complete information case) and \cite{CoEtAl08} (complete uncertainty case), where we suppose to have partial information about the network to be jammed. In particular, we suppose to have an estimate of the SIR balance in each testpoint of the network and we propose an original Robust Optimization (RO) approach to provide protection against estimated deviations. Our RO approach also presents a different way of dealing with uncertainty w.r.t. the scenario-based approach of \cite{CoEtAl09};
\item proposing an original \emph{robust cutting-plane algorithm} to tackle the right-hand-side (RHS) uncertainty coming from uncertain SIR quantities. Tackling RHS uncertainty by a canonical row-wise uncertainty approach and cardinality-constrained uncertainty sets like \cite{BeSi04} leads to trivial and conservative robust counterparts. Our new algorithm allows to overcome these conservatism and model rigidity and to exploit in an innovative way the potential of recent  Multiband Robust Optimization (see e.g., \cite{BaEtAl13,BuDA12a,BuDA12b}).
\end{enumerate}

\smallskip

\noindent
The remainder of the paper is organized as follows. In Section \ref{sec:WND}, we introduce fundamentals of \emph{realistic} wireless network design. These concepts are then used in Section \ref{sec:NJP} to derive an optimization model for the NJP. In Section \ref{sec:Robust-NJP}, we discuss data uncertainty in jamming and present our original algorithm. Finally, in Section \ref{sec:experiments} we evaluate our original algorithm on realistic wireless network instances, to then derive conclusions in Section \ref{sec:conclusions}.

\section{Classical Wireless Network Design}
\label{sec:WND}

To define our new model for network jamming, we first discuss closely related concepts from wireless network design. For modeling purposes, a wireless network can be essentially described as a set $S$ of
\emph{transceivers stations (TRXs)} that provide a telecommunication service to a set of users that are located in a target area. In line with recommendations by telecommunication regulatory bodies (e.g., \cite{AGCOM,Chester}), we decompose the target area into a set $T$ of \emph{testpoints (TPs)}, namely elementary portions of territory of identical and squared size. Each TP is assumed to correspond to a ``superuser'' that is representative of all users within the corresponding elementary area.

TRXs and TPs are characterized by a location (geographical coordinates) and a number of radio-electrical
parameters (e.g., power emission, frequency channel, transmission scheme).  The \emph{Wireless
Network Design Problem}  (WND)
consists in establishing the location and suitable values for the parameters of
the TRXs to optimize
an objective function that expresses the interest of the
decision maker (e.g., maximizing a service revenue function).
For an exhaustive introduction to the WND, we refer the reader to
\cite{DA12,DAMaSa13,KeOlRa10,RePa06}.

An optimization model for the WND typically focuses attention only on a subset of the parameters characterizing a TRX. In particular, the majority of the models considers the setting of power emission of TRXs and the assignment of served TPs to TRXs as the main decision variables. These are indeed two critical decisions that must be taken by a network administrator, as indicated in several real studies (e.g., \cite{CaEtAl11,DA12,DAMaSa11,DAMaSa13,KeOlRa10,MaRoSm06,RePa06}). Other parameters that are commonly considered are the frequency and the transmission scheme used to serve a terminal (e.g., \cite{DAMa09,DAMaSa11,MoSmAl04}). In \cite{DA12,MaRoSm06}, several distinct versions of the WND are presented and a hierarchy of WND problems is identified.

Let us now focus attention on a TP $t \in T$: when covered with service, $t$ is served by a single TRX $s \in S$, called \emph{server}, that provides the telecommunication service to it. Once the server of a TP is chosen, all the other TRXs are \emph{interferers} and reduce the quality of service obtained by $t$ from its server $s$.
Analytically, if we denote by $p_s > 0$ the power emission of a TRX $s \in S$, a TP $t \in T$ is covered with service (or
\emph{served}) when the ratio of the \emph{received} service power to the sum of the \emph{received} interfering powers
(\emph{signal-to-interference ratio} - \emph{SIR}) is above a threshold $\delta >
0$, which depends on the desired quality of service \cite{Ra01}:
\begin{equation}
\label{eq:firstSIRineq}
SIR_{t s}(p) = \frac{a_{t s(t)} \cdot p_{s(t)}}
{N + \sum_{s \in S\setminus\{s(t)\}} a_{ts} \cdot p_s}
\hspace{0.1cm}
\geq
\hspace{0.1cm}
\delta \;  .
\end{equation}

\noindent
In this inequality:
1) $s(t) \in S$ is the server of TP $t$;
2) $N > 0$ represents the noise of the system, which is commonly regarded as a constant whose value depends upon the frequency used for transmissions (see \cite{DAMaSa13,Ra01});
3) the \emph{received  power} $P_{s}(t)$ that $t$ gets from any TRX $s \in S$ is the product of the power $p_s$ emitted by $s$ multiplied by the factor $a_{ts}$, i.e.
$P_{s}(t) = a_{ts}\cdot p_s$.
The factor $a_{ts}$ is called {\em fading
coefficient}, lies in the range $[0,1]$ and summarizes the reduction in power that a signal experiences while propagating from $s$ to $t$ \cite{Ra01}.

Through simple algebra operations, inequality (\ref{eq:firstSIRineq}) can be
transformed into the
following linear inequality, commonly called \emph{SIR inequality}:
\begin{equation}\label{eq:secondSIRineq}
a_{ts(t)}  \cdot p_{s(t)} - \delta  \sum_{s \in S \setminus\{s(t)\}} a_{t s}  \cdot p_s
\hspace{0.1cm}
\geq
\hspace{0.1cm}
\delta \cdot N \;  .
\end{equation}

\noindent
Since assessing service coverage is a central question in the design of any
wireless network, the SIR inequality constitutes the core of any optimization
problem used in wireless network design.
In a hierarchy of WND problems,
a particularly relevant case is constituted by the \emph{Scheduling and Power Assignment Problem (SPAP)} \cite{DA12,DAMaSa11,DAMaSa13,MaRoSm06,RePa06}.
In the SPAP, two decisions must be taken: 1) setting the power emission of each TRX $s \in S$ and 2) assigning served TPs to activated TRXs (note that this corresponds to identify a subset of service links TRX-TP that can be scheduled simultaneously without interference, so we use the term \emph{scheduling}).
To model these two decisions, two types of decision variables are commonly introduced:
\begin{enumerate}
   \item a non-negative continuous \emph{power variable} $p_s \in [0,P_{TRX}]$ to represent the power emission of each TRX $s \in S$;
   \item a binary \emph{service assignment variable} $x_{ts} \in \{0,1\}$, $\forall \hspace{0.1cm} t \in T, s \in S$, that is equal to 1 if TP $t \in T$ is served by TRX $s \in S$ and 0 otherwise.
\end{enumerate}

\noindent
By exploiting these two families of variables, the SPAP can be naturally formulated as the following Mixed-Integer Linear Program (SPAP-MILP):
\begin{align}
\max\;\;
&
\sum_{t \in T} \sum_{s \in S} r_{t} \cdot x_{ts}
&&
\tag{SPAP-MILP}
\\
&
a_{ts}  \cdot p_{s} - \delta  \sum_{s \in S \setminus\{\sigma\}} a_{t \sigma}  \cdot p_\sigma
+ M  \cdot  (1 - x_{ts})
\hspace{0.1cm}
\geq
\hspace{0.1cm}
\delta \cdot N
&&
t \in T, s \in S
\label{SPAPsir}
\\
&
\sum_{s \in S} x_{ts} \leq 1
&&
t \in T
\label{SPAPgub}
\\
&
0 \leq p_s \leq P_{TRX}
&&
s \in S
\nonumber
\\
&
x_{ts} \in \{0,1\}
&&
t \in T, s \in S \; .
\nonumber
\end{align}

\noindent
The objective function aims at maximizing the total revenue obtained by serving testpoints (the coverage of each TP generates a revenue equal to $r_t > 0$). Each constraint \eqref{SPAPsir} corresponds with the SIR coverage condition \eqref{eq:firstSIRineq} defined for a TP $t$ served by TRX $s$ and includes a sufficiently large value $M$ (so-called, \emph{big-M coefficient}) to activate/deactivate the constraint.
Finally, constraints \eqref{SPAPgub} impose that each TP is served by at most one TRX.

\section{The Wireless Network Jamming Problem}
\label{sec:NJP}

Consider a wireless network designed by solving SPAP-MILP.
Our aim is now to compromise the functionality of the network by deploying jamming stations (\emph{jammers}). A jammer has the essential task of emitting a signal on the same frequency channel used by the jammed network to interfere with the transmissions of the TRXs and destroy their service.

Let $J$ be the set of deployed jammers and denote by $p_j > 0$ the power emission of each jammer $j \in J$. The presence of the jammers in the wireless network has the effect of creating an \emph{additional interfering summation} in the SIR inequality \eqref{eq:secondSIRineq} associated with each testpoint $t \in T$, namely:
\begin{equation}
\label{eq:secondjammedSIRineq}
a_{ts(t)}  \cdot p_{s(t)}
\hspace{0.1cm} - \delta \sum_{s \in S \setminus\{s(t)\}} a_{t s} \cdot p_s
\hspace{0.1cm} - \delta \sum_{j \in J} a_{t j} \cdot p_j
\hspace{0.1cm}
\geq
\hspace{0.1cm}
\delta \cdot N \;  .
\end{equation}

\noindent
Assume now that we want to interdict the communications in the network by jamming.
To operate the jamming we are allowed to choose the subset $J' \subseteq J$ of jammers that are activated and the corresponding power emissions $p_j \in [0,P_{JAM}]$, $\forall j \in J'$. We stress that it is rational to set the power emission of each activated jammer to its highest feasible value, since this provides the highest jamming effect. So we assume that if $j \in J$ is activated, then it emits at maximum power, i.e. $p_j = P_{JAM}$.

Let us consider now a wireless network made up of a set of TRXs $S$ providing the service to a set of TPs $T$. Moreover, let us assume that this network has been configured by solving problem SPAP-MILP. So we have at disposal a feasible solution $(\bar{x},\bar{p})$ of SPAP-MILP, which identifies the subset $T'\subseteq T$ of served TPs (i.e., $T' = \{t \in T: \bar{x}_{ts} = 1 \mbox{ for some } s \in S\}$ and the power emission $\bar{p}_{s}$ of each TRX $s \in S$.
Given a served TP $t \in T'$, we know that the corresponding SIR inequality \eqref{eq:secondjammedSIRineq} without the jamming terms is satisfied by the feasible power vector $\bar{p}$, i.e.
$
a_{ts(t)}  \cdot \bar{p}_{s(t)}- \delta \sum_{s \in S \setminus\{s(t)\}} a_{t s} \cdot \bar{p}_s
\geq
\delta \cdot N \; .
$

\noindent
In order to compromise service in $t$ by jamming, the SIR inequality must be violated and we must activate a subset $J' \subseteq J$ of jammers such that:
\begin{equation}
\label{eq:violatedSIRineq}
a_{ts(t)}  \cdot \bar{p}_{s(t)}
\hspace{0.1cm} - \delta \sum_{s \in S \setminus\{s(t)\}} a_{t s} \cdot \bar{p}_s
\hspace{0.1cm} - \delta \sum_{j \in J'} a_{t j} \cdot P_{JAM}
\hspace{0.1cm}
<
\hspace{0.1cm}
\delta \cdot N \;  .
\end{equation}

\noindent
If we introduce a binary \emph{jammer activation variable} $y_{j} \in \{0,1\}$ for each $j \in J$, which is equal to 1 if jammer $j$ is activated (at its maximum power $P_{JAM}$, as discussed above) and 0 otherwise, and if we define the quantity:
\begin{equation}
\label{eq:SIRbalance}
\Delta SIR_{ts(t)} (\bar{p})
=
a_{ts(t)}  \cdot \bar{p}_{s(t)}
- \delta \sum_{s \in S \setminus\{s(t)\}} a_{t s} \cdot \bar{p}_s
- \delta \cdot N \;  ,
\end{equation}

\noindent
which expresses the SIR balance in TP $t$ when assigned to TRX $s(t)$ for a power vector $\bar{p}$, then the violated SIR inequality \eqref{eq:violatedSIRineq} can be rewritten as:
\begin{equation}
\label{eq:violatedSIRineq_concise}
\delta \sum_{j \in J} a_{t j} \cdot P_{JAM} \cdot y_{j}
\hspace{0.1cm}
>
\hspace{0.1cm}
\Delta SIR_{ts(t)} (\bar{p}) \;  .
\end{equation}

\noindent
This inequality expresses the jamming condition: to jam and deny service in a TP, we must activate a subset of jammers whose total power received in the TP is greater than the SIR balance granted by the TRXs of the wireless network.
This inequality constitutes the central element of the new jamming optimization model that we introduce in the next paragraph.

\subsection{A SIR-based model for the Wireless Network Jamming Problem}

In our study, given an operating wireless network, we define the \emph{Wireless Network Jamming Problem} (NJP) as follows:
we must select which jammers to activate to maximize a profit function associated with jamming of served TPs, while respecting a budget that we have at disposal for the activation. The budget is introduced to model the fact that in real-world deployments we expect to have limited resources available, thus restricting the possibility of deploying jammers in the target area.

Suppose now that for each potentially activable jammer $j \in J$, we have the possibility of choosing among $m \in M = \{1,\ldots,|M|\}$ typologies of jamming devices, each associated with a distinct maximum power emission $P_{JAM}^{m}$ and a distinct cost of deployment $c_{j}^{m} > 0$. In particular, $\forall j \in J$ we assume that $P_{JAM}^{m} < P_{JAM}^{m+1}$ and $c_{j}^{m} < c_{j}^{m+1}$, $\forall m \in \{1,\ldots,|M|-1\}$.

If we add an index $m \in M$ to the \emph{jammer activation variables} to consider the presence of multiple jamming devices and we introduce binary \emph{jamming variables}  $z_{t} \in \{0,1\}$, $\forall \hspace{0.1cm} t \in T$, that are equal to 1 if served TP $t \in T'$ is jammed and 0 otherwise, the NJP can be modeled as the following 0-1 linear program:
\begin{align}
\max
&
\sum_{t \in T'} \pi_{t} \cdot z_{t}
&&
\tag{NJP-01}
\nonumber
\\
&
\delta \sum_{j \in J} \sum_{m \in M} a_{t j} \cdot P_{JAM}^{m} \cdot y_{jm}
+ M \cdot (1 - z_t)
\hspace{0.05cm}
\geq
\hspace{0.05cm}
\Delta SIR_{t} + \epsilon
&&
t \in T'
\label{GenNJP_SIR}
\\
&
\sum_{j \in J} \sum_{m \in M} c_{jm} \cdot y_{jm} \leq C
&&
\label{GenNJP_costbudget}
\\
&
\sum_{m \in M} y_{jm} \leq 1
&&
j \in J
\label{GenNJP_GUB}
\\
&
z_{t} \in \{0,1\}
&&
t \in T'
\nonumber
\\
&
y_{jm} \in \{0,1\}
&&
j \in J, m \in M
\nonumber
\end{align}

\noindent
In this model, we maximize an objective function that includes profits $\pi_{t} > 0$ deriving from jamming served TPs $t \in T'$. Constraints \eqref{GenNJP_SIR} are derived from the SIR jamming condition \eqref{eq:violatedSIRineq_concise}. Note that the constraints include a big-M term for activation/deactivation: this is necessary since, due to the budget constraint \eqref{GenNJP_costbudget}, it may happen that not all served TPs can be jammed at the same time, thus requiring to choose those that are jammed.
Constraint \eqref{GenNJP_costbudget} expresses the budget condition: we can activate a subset of jammers whose total cost does not exceed the total budget $C>0$. Finally, constraints \eqref{GenNJP_GUB} impose that we can install at most one jamming device in each activated jammer.

\smallskip
\noindent
\textbf{Remark.} In constraints \eqref{GenNJP_SIR}, we just show the dependence of the SIR balance $\Delta SIR_{t}$ on the TP index t,
omitting
$s(t)$ and $\bar{p}$.
We do this since, assuming the point of view of the NJP decision maker, we are only interested in knowing the value of the SIR balance in $t$ and we can neglect the information about the serving TRX and the power of the TRXs. We also highlight the presence of a very small value $\epsilon > 0$ to overcome the strict inequality of \eqref{eq:violatedSIRineq_concise}.

\section{Multiband Robust Optimization in Wireless Network Jamming}
\label{sec:Robust-NJP}

In the previous section, we have considered a \emph{deterministic} version of the NJP, namely we have assumed to know exactly the value of all data involved in the problem. However, in practice this assumption is likely to be not true, as also discussed in \cite{CoEtAl08,CoEtAl09}:
assuming to possess a complete knowledge about the unfriendly network is unrealistic, especially in defence and security applications, where it may be very difficult or even dangerous to gather accurate information. Instead it is rational to assume that we can just rely on estimates of the position and the radio-electrical configuration of the TRXs.
As a consequence, it is highly reasonable to assume that we just possess an estimate of the value of the SIR balance $\Delta SIR_{t}$ in every TP $t$.
Following a practice that we have observed among wireless network design professionals dealing with uncertain SIR quantities (see \cite{DA12}), we use an estimate $\Delta\overline{SIR}_{t}$ as a reference \emph{nominal value} to define an interval of variation of the quantity, whose bounds reflect the reliability of the limited information that we have at disposal and our risk aversion.
If we denote by $d_{t}^{-} < 0$ and $d_{t}^{+} > 0$ the maximum negative and positive deviation from $\Delta\overline{SIR}_{t}$ that we have derived on the basis of our limited information, then the (unknown) actual value $\Delta SIR_{t}$ belongs to the interval:
$
[\Delta\overline{SIR}_{t} + d_{t}^{-}, \hspace{0.2cm} \Delta\overline{SIR}_{t} + d_{t}^{+}] \;.
$
We note that the definition of the interval of variation of $\Delta SIR_{t}$ can also take into account the intrinsic uncertainty of the fading coefficients $a_{t j}$ of the jammers: propagation of wireless signals in a real environment is affected by many distinct factors (e.g., distance between the TRX and the TP, presence of obstacles, weather) that are very hard to assess precisely. Therefore the exact value of the fading coefficients is typically not known
(see \cite{DA12} and \cite{Ra01} for an exhaustive discussion).

An example may help to clarify the negative effects of uncertainty in the NJP. Note that, as common in the WND practice, we express fading and power quantities according to a \emph{decibel} (\emph{dB}) scale. More specifically, since we measure power quantities in \emph{milliwatts} (\emph{mW}), we express power in decibels by referring to \emph{decibel-milliwatts} (\emph{dBmW}).

\noindent
\textbf{Example 1 (Uncertainty in the NJP).}
Consider a TP that is part of a wireless network subject to a noise of $N_{dB} = -114$ $dBmW$ and operating with a SIR threshold $\delta_{dB} = 10$ $dB$ and that receives a serving power of $-48$ $dBmW$ and a total interfering power (including noise) of $-61$ $dBmW$. By formula \eqref{eq:firstSIRineq}, the SIR in the TP is higher than $\delta_{dB}$, being equal to about $13$ $dB$. Therefore the TP is served.
The corresponding SIR balance $\Delta SIR$ can be computed by formula \eqref{eq:SIRbalance} and is equal to about $-50$ $dBmW$.

Suppose now that we want to jam the TP and that we can install a single jammer in a site associated with a fading coefficient of $-77$ $dB$ towards the considered TP. Additionally, suppose that the jammer can accommodate either a device $J_1$ with $P_{JAM}^{J_1} = 20$ $dBmW$ or a more powerful and costly device $J_2$ with $P_{JAM}^{J_2} = 27$ $dBmW$.
If we assume to know all the features of the jammed network, then we can successfully jam the TP by generating an \emph{additional} received interfering power of at least about $17$ $dBmW$.
So installing the less powerful device $J_1$ emitting at $P_{JAM}$ is sufficient to deny service.

However, as previously discussed, in real-world scenarios it is likely that we do not know the exact value of $\Delta SIR$, but we just possess an estimate $\Delta\overline{SIR}$ and an interval of deviation. Suppose then that our estimate is $\Delta\overline{SIR} = - 50$ $dBmW$ and that we consider reasonable to experience deviations up to $\pm$20\% of this value. So the actual value $\Delta SIR$ lies in the interval $[-60,-40]$ dBmW. This interval reflects how trustable we consider the available information about the unfriendly network and expresses also our personal risk aversion.
If the worst deviation occurs, we have $\Delta SIR = - 40$ $dBmW$ and activating $J_1$ would be no more sufficient to successfully jam the TP: the jamming solution deploying $J_1$ would be infeasible and thus useless. So, if we want to be protected against this deviation, we should switch to the more powerful jammer $J_2$, at the price of a higher deploying cost.
\qed

\medskip
\noindent
As the example has shown, the presence of uncertain data in an optimization problem can lead to very bad consequences: as it is known from sensitivity analysis, even small deviations in the value of input data may make an optimal solution heavily suboptimal, whereas feasible solutions may reveal to be infeasible and thus completely useless in practice \cite{BeElNe09,BeBrCa11}. In our application, it is thus not possible to optimize just referring to the nominal values $\Delta\overline{SIR}_{t}$, but we must take into account the possibility of deviations in an interval.

Many methodologies have been proposed over the years to deal with data uncertainty: \emph{Stochastic Programming} is commonly considered the oldest, while in the last decade \emph{Robust Optimization} has known a wide success, especially in real-world applications thanks to its accessibility and computational tractability. We refer the reader to \cite{BeElNe09,BeBrCa11} for a general discussion about the impact of data uncertainty in optimization and for an overview of the main methodologies proposed in literature to deal with uncertain data. The two references are in particular focused on theory and applications of Robust Optimization (RO), the methodology that we adopt in this paper to tackle data uncertainty. RO is based on two main concepts: 1) the decision maker defines an \emph{uncertainty set}, which reflects his risk aversion and identifies the deviations of coefficients against which he wants to be protected; 2) protection against deviations specified by the uncertainty set is guaranteed under the form of hard constraints that cut off all the feasible solutions that may become infeasible for some deviations included the uncertainty set.
Formally, suppose that we are given a generic 0-1 linear program:
$$
v
\hspace{0.1cm} = \hspace{0.1cm}
\max
\hspace{0.1cm}
c' \hspace{0.05cm} x
\hspace{0.5cm}
\mbox{ with }
\hspace{0.1cm}
x \in {\cal F}
=
        \{
        A \hspace{0.05cm} x \leq b,
        \hspace{0.15cm}
        x \in \{0,1\}^{n}
        \}
$$

\noindent
and that the coefficient matrix $A$ is uncertain, i.e. we do not know the exact value of its entries. However, we are able to identify a family $\cal A$ of coefficient matrices that represent possible valorizations of the uncertain matrix $A$, i.e. $A \in \cal A$. This family represents the uncertainty set of the robust problem. A \emph{robust optimal solution}, i.e. a solution protected against data deviations, can be computed by solving the \emph{robust counterpart} of the original problem:
$$
v^{{\cal R}}
\hspace{0.1cm} = \hspace{0.1cm}
\max
\hspace{0.1cm}
c' \hspace{0.05cm} x
\hspace{0.5cm}
\mbox{ with }
\hspace{0.1cm}
x \in {\cal R}
=
        \{
        A \hspace{0.1cm} x \leq b
        \hspace{0.3cm}
        \forall A \in {\cal A},
        \hspace{0.15cm}
        x \in \{0,1\}^{n}
        \} \; ,
$$

\noindent
which considers only the solutions that are feasible for all the coefficient matrices in the uncertainty set ${\cal A}$. Therefore, the robust feasible set is such that ${\cal R} \subseteq {\cal F}$. The choice of the coefficient matrices included in ${\cal A}$ should reflect the risk aversion of the decision maker.

Guaranteeing protection against data deviations entails the so-called \emph{price of robustness} \cite{BeSi04}: the optimal value of the robust counterpart is in general worse than the optimal value of the original problem, i.e., $v^{{\cal R}} \leq v$, due to having restricted the feasible set to only robust solutions. The price of robustness reflects the features of the uncertainty set: uncertainty sets expressing higher risk aversion will take into account more severe and unlikely deviations, leading to higher protection but also higher price of robustness; conversely, uncertainty sets expressing risky attitudes will tend to neglect improbable deviations, offering less protection but also a reduced price of robustness.
\\
In the next paragraph, we fully describe the uncertainty model that we adopt.

\subsection{RHS Uncertainty in Wireless Network Jamming}

The data uncertainty affecting our  problem needs a special discussion. As pointed out in \cite{BeElNe09,BeBrCa11}, most RO models considers so-called \emph{row-wise} uncertainty. This means that protection against data deviations is separately defined for each constraint subject to uncertainty, by considering the worst total deviation that the constraint may experience w.r.t. the uncertainty set. More formally, consider again a generic uncertain 0-1 linear program:
\begin{align*}
\max
\sum_{j\in J} c_j \cdot x_j
\hspace{0.4cm}
\mbox{ s.t. }
\sum_{j\in J} a_{ij} \cdot x_j \leq b_i
\hspace{0.3cm}
i\in I,
\hspace{0.4cm}
x_j \in \{0,1\}
\hspace{0.3cm}
j \in J
\; .
\end{align*}

\noindent
where w.l.o.g we assume that the uncertainty just regards the coefficients $a_{ij}$ (uncertainty affecting cost coefficients or RHSs can be easily reformulated as coefficient matrix uncertainty, see \cite{BeSi04}).
If we denote the uncertainty set by $U$, following a row-wise uncertainty paradigm the robust counterpart is:
\begin{align*}
\max
\sum_{j\in J} c_j \cdot x_j
\hspace{0.15cm}
\mbox{ s.t. }
\sum_{j\in J} a_{ij} \cdot x_j + DEV_i(x,U) \leq b_i
\hspace{0.25cm}
i\in I,
\hspace{0.30cm}
x_j \in \{0,1\}
\hspace{0.15cm}
j \in J
\; .
\end{align*}

\noindent
where each uncertain constraint $i \in I$ 1) refers to the nominal value $\bar{a}_{ij}$ of each coefficient and 2) includes an additional term $DEV_i(x,U)$ to represent the maximum total deviation that $i$ may experience for the solution $x$ and the uncertainty set $U$. This problem is actually non-linear since $DEV_i(x,U)$ hides a maximization problem based on the uncertainty set definition (see \cite{BeSi04,BuDA12a,BuDA12b}).

A central question in RO is how to model the uncertainty through a suitable uncertainty set $U$. The majority of applied studies of RO model $U$ as a cardinality-constrained uncertainty set  \cite{BeBrCa11}, primarily referring to the renowned $\Gamma$-robustness model ($\Gamma$-Rob) by Bertsimas and Sim \cite{BeSi04}. The main feature of these particular uncertainty sets is to impose an upper bound on the number of coefficients that may deviate to their worst value in each constraint. The non-linearity of the robust counterpart due to the presence of $DEV_i(x,U)$ is then solved by exploiting strong duality and defining a larger but compact and linear robust counterpart, as explained in \cite{BeSi04} and \cite{BuDA12a,BuDA12b}.

In relation to this general row-wise RO setting, the uncertain NJP that we consider is a special type of uncertain problem: uncertainty just affects the RHS of each SIR constraint of NJP-01. As a consequence, if we adopt row-wise uncertainty and a cardinality-constrained uncertainty set, then the upper bound on the number of deviating coefficient in each constraint \eqref{GenNJP_SIR} is equal to either 0 or 1. In other words, either the constraint is not subject to uncertainty and thus the actual value and the nominal value coincide (i.e., $\Delta SIR_{t} =\Delta\overline{SIR}_{t}$) or the constraint is subject to uncertainty and thus the actual value is equal to the highest deviating value (i.e., $\Delta SIR_{t} = \Delta\overline{SIR}_{t} + d_{t}^{+}$). Thus the robust counterpart simply reduces to a nominal problem with modified RHS values. We stress that this is a very rigid representation of the uncertainty and we would like to benefit from a richer representation.

A source of inspiration for a richer model can be represented by \emph{Multiband Robust Optimization} (MB) and related multiband uncertainty sets, introduced by B\"using and D'Andreagiovanni in 2012 to generalize and refine classical $\Gamma$-Rob (see e.g., \cite{BuDA12a,BuDA12b,BuDA14} and \cite{BaEtAl13}).
In our case, we want to adopt a distinct but similar definition of multiband uncertainty.
To define this multiband-like uncertainty set for RHS uncertainty:
\begin{enumerate}
\item  we partition the overall deviation range
    $
    [d_{t}^{-}, \hspace{0.1cm} d_{t}^{+}]
    $
    into $K$ bands, defined on the basis of $K$ deviation values:
      \\
        $
        -\infty<
        {d_{t}^{-} = d_{t}^{K^{-}}<\cdots<d_{t}^{-1}
        \hspace{0.1cm}<\hspace{0.2cm}d_{t}^{0}=0\hspace{0.2cm}<\hspace{0.1cm}
        d_{t}^{1}<\cdots<d_{t}^{K^{+}}} = d_{t}^{+}
        <+\infty ;
        $
  \item through these deviation values, $K$ deviation bands are defined, namely:
    a set of positive deviation bands $k\in \{1,\ldots,K^{+}\}$ and  a set of negative deviation bands $k \in \{K^{-}+1,\ldots,-1,0\}$, such that a band $k\in \{K^{-}+1,\ldots,K^{+}\}$ corresponds to the range $(d_{t}^{k-1},d_{t}^{k}]$, and band $k = K^{-}$
    corresponds to the single value $d_{t}^{K^{-}}$. Note that $K = K^{+} + K^{-}$;
  \item considering the RHS values $\Delta SIR_{t}$ of the entire set of constraint \eqref{GenNJP_SIR}, we impose a lower and upper bound on the number of values that may deviate in each band: for each band $k\in K$, we introduce two bounds $l_k, u_{k} \in \mathbb{Z}_{+}$: $0 \leq l_k \leq u_{k} \leq |T'|$. As additional assumptions, we do not limit the number of coefficients that may deviate in band $k = 0$ (i.e., $u_{0}=|T'|$), and we impose that $\sum_{k \in K} l_k \leq |T'|$, to ensure the existence of a feasible realization of the uncertainty set.
\end{enumerate}

\noindent
We call this set \emph{RHS-Multiband Set (RHS-MB)}.

\smallskip
\noindent
\textbf{Remark.}
We stress that point 3 differs from the classical definition of multiband uncertainty set, presented in  \cite{BuDA12a,BuDA12b}, where a row-wise uncertainty perspective is assumed and the system of bounds for the bands is defined separately for each uncertain constraint of the problem.

\smallskip
\noindent
An MB uncertainty set is particularly suitable to represent the past behaviour of uncertainty represented by histograms, as explained in \cite{BaEtAl13,BuDA12a,BuDA12b}. Moreover, such set has the advantage of taking into account negative deviation bands, which are neglected in classical cardinality-constrained sets: we want of course to be protected against positive deviations that lead to infeasibility, but in real-world applications we commonly experience also negative deviations,
which compensate the positive deviations and reduce the conservatism of solutions.

A critical question is now: how can we solve the uncertain NJP when RHS uncertainty is modeled by a RHS-Multiband Set? In the case of a purely linear program, we could define the dual problem of our uncertain problem, thus transforming the RHS uncertainty into objective function uncertainty and then adopt a standard RO dualization approach and reach a compact robust counterpart, as in \cite{Mi09}.
However, due to the integrality constraints, the classical dualization approach in our case cannot be operated.

As an alternative, we can adopt a \emph{robust cutting plane approach}: we solve NJP-01 obtaining an optimal solution, then we check whether the solution is also robust and feasible w.r.t. a specified RHS-MB. If this is the case, we have found a robust optimal solution and we have done. Otherwise we separate a \emph{robustness cut}, namely an inequality that cut off this non-robust solution, we add the cut to the problem and we solve the new resulting problem. This basic step is then iterated as in a canonical cutting-plane algorithm, until no new cut is separated and thus the generated solution is robust and optimal.

Under canonical row-wise uncertainty, in $\Gamma$-Rob and MB, robustness cuts can be efficiently separated.
For $\Gamma$-Rob, the separation of a robustness cut is trivial and just consists in sorting the deviations and choosing the worst $\Gamma > 0$ \cite{FiMo12}. This simple approach is not valid instead for MB, but in \cite{BuDA12a,BuDA12b} we proved that the separation can be done in polynomial time by solving a min-cost flow problem.

As we stressed above, RHS-MB poses a new challenge.
More formally, suppose that we have a feasible solution $(\bar{z},\bar{y})$ to NJP-01 and that we want to check its robustness. Let us denote by $T'$ the subset of TPs that are successfully jammed.
A robustness cut is generated by solving the following 0-1 linear program, that can be interpreted as the problem of an adversarial that attempts to compromise the feasibility of our optimal jamming solution by picking up the worst deviation for $(\bar{z},\bar{y})$ allowed by RHS-MB.
\begin{align}
V
=
\max
&
\sum_{t \in T'} v_t
&&
\tag{SEP}
\label{MBdev_objFunction}
\\
&
\sum_{k\in K} d_{t}^{k} \cdot w_{t}^{k}
+ M \cdot (1 - v_t) \geq JAM_t - \Delta\overline{SIR}_{t}
\label{MBdev_JAMviolation}
&&
t \in T'
\\
&
l_{k} \leq \sum_{j\in J} w_{t}^{k} \leq u_{k}
&&
k\in K
\label{MBdev_constraint}
\\
&
\sum_{k\in K} w_{t}^{k} \leq 1
&&
t \in T'
\label{MBgub_constraint}
\\
&
w_{t}^{k} \in \{0,1\}
&&
t \in T', k \in K
\label{MBgub_variables}
\\
&
v_{t} \in \{0,1\}
&&
t \in T'
\label{MBgub_deniedJamming}
\; .
\end{align}

\noindent
The separation problem SEP uses 1) a binary variable $v_{t}$, $\forall t \in T'$ that is equal to 1 when the jamming of TP $t$ is denied and 0 otherwise; 2) a binary variable $w_{t}^{k}$ that is equal to 1 when the SIR balance $\Delta SIR_{t}$ of $t$ deviates in band $k$ and 0 otherwise.
The objective function aims at maximizing the number of TPs whose jamming is denied by the adversarial. A constraint \eqref{MBdev_JAMviolation} expresses the violation of the corresponding constraint \eqref{GenNJP_SIR} when the jamming of TP $t$ is denied by feasible deviations of the SIR balance according to RHS-MB, namely
$JAM_t < \Delta SIR_{t} + \sum_{k\in K} d_{t}^{k} \cdot w_{t}^{k}$, where $JAM_t = \delta \sum_{j \in J} \sum_{m \in M} a_{t j} \cdot P_{JAM}^{m} \cdot \bar{y}_{jm}$ is the total jamming power that $t$ receives for jamming solution $(\bar{z},\bar{y})$.
Constraints (\ref{MBdev_constraint})-(\ref{MBgub_constraint}) enforce the structure of the uncertainty set RHS-MB: the first family imposes the lower and upper bounds on the number of RHS values $\Delta SIR_{t}$ that may deviate in each band $k \in K$, whereas the second family imposes that each value $\Delta SIR_{t}$ deviates in at most one band (note that $\sum_{k\in K} w_{t}^{k} = 0$ corresponds with no deviation and is equivalent to $w_{t}^{0} = 1$).

It is easy to observe that if the optimal value $V$ of SEP is equal to 0, then $(\bar{z},\bar{y})$ is robust, since it is not possible to compromise the jamming of any TP for the given uncertainty set RHS-MB. On the contrary, if $V \geq 1$ and $(v^{*},w^{*})$ is an optimal solution of SEP, then $(\bar{z},\bar{y})$ is not robust, the jamming of $V$ TPs may be compromised and
\begin{equation}
\sum_{t \in T': \hspace{0.1cm} v_t^{*}=1} z_t \leq V - 1
\label{eq:robCut}
\end{equation}
is evidently a robustness cut that we must add to the original problem. After this we can iterate the basic robustness check step.

The general structure of the proposed robust cutting plane algorithm is described in Algorithm 1.
Assuming to use a solver like CPLEX implementing a branch-and-cut solution algorithm, the separation problem is solved every time that the solver finds a feasible solution to NJP-01. If a robustness cut is identified for the current solution, then it is added as constraint to NJP-01.
\begin{algorithm}
\caption{Robust Cutting Planes for NJP subject to RHS-MB uncertainty}
\label{ALGhybrid}
\begin{algorithmic}[1]
\Require an instance of NJP-01 and of RHS-MB
\Ensure a robust optimal solution $(z^{*},y^{*})$ to NJP-01 w.r.t. RHS-MB (if existent)
\State Solve NJP-01 by a branch-and-cut-based MIP solver (denoted by SOLVER)
\While{SOLVER has not find a robust optimal solution $(z^{*},y^{*})$ to NJP-01 or has proved that $(z^{*},y^{*})$ does not exist}
    \State Run SOLVER
    \If{SOLVER finds a feasible solution $(\bar{z},\bar{y})$ to NJP-01}
        \State Solve SEP for $(\bar{z},\bar{y})$ and RHS-MB
        \If{$V > 0$}
            \State Generate a robustness cut \eqref{eq:robCut} and add it to NJP-01
        \EndIf
    \EndIf
\EndWhile
\end{algorithmic}
\end{algorithm}

\section{Computational results}
\label{sec:experiments}

To evaluate the performance of our original robust cutting plane algorithm, we considered a set of 15 realistic instances, based on network data defined in collaboration with network engineers of the Technical Strategy \& Innovations Unit of British Telecom Italia. All the instances refer to a fixed WiMAX network (see \cite{AnGhMu07}, \cite{DAMa09} and \cite{DA12} for an introduction to WiMAX technology and modeling) and are based on real terrain data model and population statistics of a residential urban area from the administrative district of Rome (Italy). The instances consider distinct networks with up to $|T| = 224$ TPs and $|S| = 20$ TRXs, operating on one of the frequency channels reserved for WiMAX transmissions in Italy in the band [3.4$\div$3.6] GHz and using a QAM-16 modulation scheme.
We used these data to build the MILP problem SPAP-MILP for each instance and obtain realistic wireless network configurations to jam by solving the uncertain version of problem NJP-01. The revenue $r_t$ associated with the service coverage of each TP was derived from population statistics.

In order to build NJP-01 and set the robust cutting-plane algorithm, we assume that we know the set $T'$ of served TPs. However, we also assume that we do not exactly the value of the SIR balance granted by the solution of SPAP-MILP, but we just have at disposal an estimate $\Delta\overline{SIR}_{t}$ (different from the actual value provided by the solution). On the basis of discussions on the topic with network professionals,
we decided to model deviations through an RHS multiband uncertainty set including 5 deviation bands (2 negative and 2 positive, besides the null deviation band) and with a basic deviation of each band equal to 20\% of the nominal value. Concerning the jammers, we supposed to have three typologies of jammers (i.e., $|M| = 3$) with a cost of deployment reflecting the population in the TPs and increasing as the population in the TP increases (we assume a higher risk of deployment in more populated areas where the jammers could be discovered). The profit $\pi_t$ of successfully jamming a TP was also based on population data.

All experiments were made on a 2.70 GHz Intel Core i7 with 8 GB. The code was written in the C/C++ programming language and the optimization problems were solved by IBM ILOG CPLEX 12.5 with the support of Concert Technology. The results of the experiments are reported in Table \ref{tab:1},
in which the first column states the instance ID, whereas the following four columns report: the number $|T|$ of TPs and the number $|S|$ of TRXs in the SPAP-MILP instance; the number  $|T^{*}|$ of covered TPs in the feasible solution of SPAP-MILP used for building the corresponding NJP-01 instance; the number $|J|$ of jammers in the NJP-01 instance. The following four columns report instead: the optimal number \#JAM(Nom) of jammed TPs for the nominal NJP-01 problem (no uncertainty considered); the optimal number \#JAM(Rob) of jammed TPs for the robust version of NJP-01 solved by Algorithm 1; the percentage price of robustness PoR\%; the number \#Cuts of robust cuts generated during the execution of Algorithm 1.

The main observations about the results are related to the comparison between the optimal value of  the nominal problem and that of its robust version. Concerning this central point, we can observe that the price of robustness that we must face keep contained, reaching an average value of -17,1\% and a peak of -23.1\% in the case of instance I9. We consider this a reasonable price to pay to obtain the protection against the deviations that the decision maker considers relevant. Furthermore, we can notice that the number of robust cuts that are separated during the execution of Algorithms 1 is limited, especially in the case of the smaller instances.
Concerning solution time, while solving the uncertain version of NJP-01 required a time ranging from about 30 to about 70 minutes, depending upon the features of the wireless network configuration to be jammed identified by a solution of SPAP-MILP, the execution time of Algorithm 1 could reach approximately 3 hours. We believe that this time could be sensibly reduced by studying a stronger separation model and more efficient separation algorithms.

\begin{table}
\caption{Experimental results}
\label{tab:1}
\begin{tabular}{c cccc cccc}
\hline
\noalign
{\smallskip}
ID & $|T|$ & $|S|$ & $|T^{*}|$ & $|J|$ & \#JAM(Nom) & \#JAM(Rob) & PoR\% & \#Cuts
\\
\noalign{\smallskip}\hline\noalign{\smallskip}
I1 & 100 & 6 & 65 & 15 &
44 & 37 & -15.90 & 29
\\
I2 & 100 & 9 & 71 & 15 &
51 & 45 & -11.76 & 41
\\
I3 & 100 & 12 & 75 & 15 &
46 & 38 & -17.93 & 37
\\
I4 & 150 & 6 & 85 & 15 &
49 & 43 & -12.24 & 32
\\
I5 & 150 & 9 & 93 & 15 &
68 & 57 & -16.17 & 31
\\
I6 & 150 & 12 & 106 & 20 &
75 & 64 & -14.66 & 35
\\
I7 & 169 & 12 & 92 & 20 &
47 & 39 & -17.02 & 49
\\
I8 & 169 & 16 & 95 & 20 &
66 & 53 & -19.69 & 58
\\
I9 & 169 & 20 & 120 & 20 &
69 & 53 & -23.18 & 75
\\
I10 & 196 & 12 & 108 & 20 &
73 & 58 & -20.54 & 68
\\
I11 & 196 & 16 & 122 & 25 &
82 & 69 & -15.85 & 54
\\
I12 & 196 & 20 & 134 & 25 &
89 & 70 & -21.34 & 92
\\
I13 & 224 & 15 & 142 & 25 &
102 & 82 & -19.60 & 87
\\
I14 & 224 & 20 & 159 & 25 &
115 & 96 & -16.52 & 101
\\
I15 & 224 & 25 & 170 & 25 &
109 & 93 & -14.67 & 103
\\
\noalign{\smallskip}\hline
\end{tabular}
\end{table}

\section{Conclusions}
\label{sec:conclusions}
We considered the \emph{Wireless Network Jamming Problem}, namely the problem of optimally placing and configuring a set of jammers in order to inderdict communications of a wireless networks.
We revisited the models proposed in the seminal works by Commander et al., better highlighting the strong connections with classical wireless network design formulations. Moreover, we addressed the uncertain nature of the problem by proposing an original robust cutting plane algorithm, inspired by \emph{Multiband Robust Optimization}, to deal with the RHS uncertainty of the problem and overcome the rigidity of canonical row-wise uncertainty approaches. As future work, we plan to investigate stronger models for the problem, tackling in particular the presence of big-M coefficients and devising more effective and efficient separation algorithms.

\end{document}